\documentclass[11pt]{amsart}
\usepackage{amssymb}
\usepackage{amsmath}

\newcommand{\mpn}{\medskip\par\noindent}
\newcommand{\pn}{\par\noindent}

\newcommand{\bpn}{\bigskip\par\noindent}
\theoremstyle{definition}

\theoremstyle{remark}

\numberwithin{equation}{section}

\begin{document}
\newcommand{\Mod}[1]{\,(\text{\mbox{\rm mod}}\;#1)}
\title[Some identities of polynomials arising from Umbral calculus ]
{ Some identities of polynomials arising from Umbral calculus  }
\author{Dae San  Kim, Taekyun KIM, Seog-Hoon  Rim}
\begin{abstract}
In this paper, we study some properties of associated sequences in
umbral calculus. From these properties, we derive new and
interesting identities of several kinds of polynomials.
 \end{abstract}

 \maketitle

\section{\bf Introduction}
We recall that the Bernoulli polynomials are defined by the generating function to be
\begin{align*}
\frac{t}{e^t-1}e^{xt}=e^{B(x)t}=\sum_{n=0}^{\infty}B_{n}(x)\frac{t^{n}}{n!}, \quad (see
[7,8]),
\end{align*}
with the usual convention about replacing $B^n(x)$ by $B_n(x)$.
\\In the special case, $x=0$, $B_{n}(0)=B_{n}$ are called the $n$-th
Bernoulli numbers $(see[1-14])$.
\\For $r \in \mathbb{Z}_{+}$, the higher order Bernoulli polynomials  are also defined by the generating function to be
\begin{align*}
\left(\frac{t}{e^t-1} \right)^re^{xt}=\underbrace{\left(\frac{t}{e^t-1}\right)...\left(\frac{t}{e^t-1}\right)}_{r-times} e^{xt}=\sum_{n=0}^{\infty}B_n^{(r)}(x)\frac{t^n}{n!}.
\end{align*}
\\In the special case, $x=0$, $B_n^{(r)}(0)=B_n^{(r)}$ are called the $n$-th Bernoulli numbers of order $r$ $({see \,[5,6]})$.
From the definition of Bernoulli numbers, we note that
\begin{align*}
B_0=1,\quad (B+1)^n-B_n=\delta_{1,n},\quad (see
[7,8,10]),
\end{align*}
where $\delta_{n,k} $ is the Kronecker's symbol.
\\As is well Known, the Euler and higher-order Euler polynomials are also defined by the generating functions as
follows:
\begin{align*}
\frac{2}{e^t+1}e^{xt}=e^{E(x)t}=\sum_{n=0}^{\infty}E_{n}(x)\frac{t^{n}}{n!},
\end{align*}
and
\begin{align*}
\left(\frac{2}{e^t+1}
\right)^re^{xt}=\underbrace{\left(\frac{2}{e^t+1}\right)...\left(\frac{2}{e^t+1}\right)}_{r-times}
e^{xt}=\sum_{n=0}^{\infty}E_n^{(r)}(x)\frac{t^n}{n!},
\end{align*}
with the usual convention about replacing
$\left(E^{(r)}(x)\right)^n$ by $E_n^{(r)}(x)$ $(see[3,4,5,6])$.
\\Let $\mathbb{C}$ be the complex number field and let $\mathcal{F}$ be the set of all formal power series in the variable $t$ over $\mathbb{C}$ with
\begin{align*}
\mathcal{F}=\{f(t)=\sum_{n=0}^{\infty}\frac{a_k}{k!}t^k ~|~a_k \in\mathbb{C}\}.
\end{align*}
Let $\mathbb{P}=\mathbb{C}[t]$ and let $\mathbb{P^*}$ be the vector
space of all linear functional on $\mathbb{P}$. Now, we use the
notation $<L ~|~ p(x)>$ to denote the action of a linear functional
$L$ on a polynomial $p(x)$ $(see[4,11])$.
 \\The formal power series
\begin{align*}
f(t)=\sum_{k=0}^{\infty}\frac{a_k}{k!}t^k \in \mathcal{F},\quad (see [4,11]),
\end{align*}
define a linear functional on $\mathbb{P}$ by setting
\begin{align*}
<f(t) ~ | ~x^n>=a_n ,\quad for\quad all\quad  n\geq 0,\quad (see [4,11]).
\end{align*}
Thus, we have
\begin{align*}
<t^k ~ | ~x^n>=n!\delta_{n,k},  \quad (see [4]).
\end{align*}
Let $f_L(t)=\sum_{k=0}^{\infty}\frac{<L~|~ x^k>}{k!}t^k $. Then, we note that $<f_L(t) ~|~ x^n> = <L~|~x^n>$ and
so as linear functionals $L=f_L(t)$. It is known in
[11] that the map $L \mapsto f_L(t)$ is a vector space isomorphism
from $\mathbb{P^*}$ on to  $\mathcal{F}$. Henceforth, $\mathcal{F}$
will denote both the algebra of formal power series in $t$ and the
vector space of all linear functionals on $\mathbb{P}$ and so an
element $f(t)$ of $\mathcal{F}$ will be thought of as both a formal
power series and a linear functional. We shall call $\mathcal{F}$
the umbral algebra. The umbral calculus is the study
of umbral algebra and modern classical umbral calculus can be
described as a systematic study of the class of Sheffer sequences $(see[11])$.
 The order $O(f(t))$ of the power series $f(t)\neq
0$ is the smallest integer $k$ for which $a_k$ does not vanish. The
series $f(t)$ has a multiplicative inverse, denoted by $f(t)^{-1}$
or $\frac{1}{f(t)}$ if and only if $O(f(t))=0$. Such a series is
called invertible series. A series $f(t)$ for which $O(f(t))=1$
is called a delta series $(see [4,11])$. Let $f(t),g(t) \in
\mathcal{F}$. Then, we see that $$<f(t)g(t) ~|~ p(x)> = <f(t) ~|~
g(t)p(x)> = <g(t) ~|~ f(t)p(x)>, \quad (see [11]). $$ In [11], we
note that for all $f(t)$ in $\mathcal{F}$
\begin{align*}
f(t)=\sum_{k=0}^{\infty}\frac{<f(t)~|~x^k>}{k!}t^k ,
\end{align*}
and for all polynomials $p(x)$
\begin{align*}
p(x)=\sum_{k=0}^{\infty}\frac{<t^k~|~p(x)>}{k!}x^k .
\end{align*}
Thus, we get
\begin{align*}
p^{(k)}(x)= \frac{d^kp(x)}{dx^k}&=\sum_{l=k}^{\infty}\frac{1}{(l-k)!}<t^l~|~p(x)>x^{l-k},
\end{align*}
and
\begin{align*}
p^{(k)}(0)= <t^k~|~p(x)> \quad and <1 ~|~ p^{(k)}(x)>=p^{(k)}(0).
\end{align*}
From this, we have $t^kp(x)=p^{(k)}(x)= \frac{d^kp(x)}{dx^k},(k\geq
0)$. It is not difficult to show that $e^{yt}p(x)=p(x+y)$
$(see[4,11])$. Let $S_n(x)$ be a polynomial with deg $S_n(x)=n$, $f(t)$  a delta series and let $g(t)$ be an invertible series.
Then there exists a unique sequence $S_n(x)$ of polynomials with $<
g(t) f(t)^k ~|~S_n(x)>=n!\delta _{n,k}, (n,k \geq 0)$. The sequence
$S_n(x)$ is called the Sheffer sequence for $(g(t),f(t))$, which is
denoted by $S_n(x)\sim (g(t),f(t))$ . If $S_n(x) \sim (1,f(t))$,
then $S_n(x)$ is called the associated sequences for $f(t)$, or
$S_n(x)$ is associated to $f(t)$. If $S_n(x) \sim (g(t),t)$, then
$S_n(x)$ is called the Appell sequence for $g(t)$ or $S_n(x)$ is
Applell for $g(t)$ $(see [4,11])$. For $p(x) \in \mathbb{P}$, we
have $<\frac{e^{yt}-1}{t} ~|~p(x)>= \int_0^{y}p(u)du$  $(see [4,11])$
\\Recently, Dere and Simsek have studied umbral calculus related to special polynomials.
In this paper, we study some properties of associated sequences in umbral algebra. From these properties, we derive new and interesting identities of several kinds of
 polynomials.

\section{\bf Some identities of polynomials arising from umbral calculus }
Let $p_n(x) \sim (1,f(t))$ and $q_n(x) \sim (1,g(t))$. Then, for
$n\geq 1$, we have
\begin{align*}\tag{1}
q_n(x)=x\left (\frac{f(t)}{g(t)}\right) ^nx^{-1}p_n(x),\quad
(see[11]).
\end{align*}
Let us take $p_n(x)=(x)_n$ and $q_n(x)=x^n$. Then we see that $(x)_n
\sim (1,e^t-1)$ and $x^n \sim (1,t)$.
\\It is easy to show that
\begin{align*}\tag{2}
\left(\frac{e^t-1}{t}
\right)^n=\underbrace{\left(\frac{e^t-1}{t}\right)...\left(\frac{e^t-1}{t}\right)
}_{n-times}=\sum_{l=0}^{\infty}\frac{n!}{(l+n)!}S_2(l+n,n)t^l,
\end{align*}
where $S_2(n,k)$ is the stirling number of the second kind.
\\For $n\geq 1$, by (1), we get
\begin{align*}\tag{3}
x^n&=x\left(\frac{e^t-1}{t} \right)^nx^{-1}(x)_n
\\&=x\left(\frac{e^t-1}{t} \right)^n(x-1)_{n-1}
\\&=x\sum_{l=0}^{\infty}\frac{n!}{(l+n)!}S_2(l+n,n)t^l(x-1)_{n-1},
\end{align*}
where $(x)_n=x(x-1)...(x-n+1)$.
\\The stirling number of the first kind is defined by
\begin{align*}\tag{4}
(x)_n=\sum_{l=0}^{n}S_1(n,l)x^l, \quad (see [6,11]).
\end{align*}
By (3) and (4), we get
\begin{align*}\tag{5}
(x+1)^n&=\sum_{l=0}^{n}\frac{(n+1)!}{(l+n+1)!}S_2(l+n+1,n+1)\sum_{m=0}^{n-l}S_1(n,l+m)x^m(l+m)_l
\\&=\sum_{m=0}^{n}\sum_{l=0}^{n-m}\frac{\binom{l+m}{l}}{\binom{l+n+1}{l}} S_2(l+n+1,n+1)S_1(n,l+m)x^m,
\end{align*}
and
\begin{align*}\tag{6}
(x+1)^m=\sum_{m=0}^{n}\binom{n}{m}x^m.
\end{align*}
Therefore , by (5) and (6), we obtain the following theorem.
\par\bigskip
{\bf  Theorem 1 }. For $m,n\in \mathbb{Z}_+=\mathbb{N} \bigcup
\{0\}$ with $ n\geq m\geq 0$ , we have
\begin{equation*}
\begin{split}
\binom{n}{m}=\sum_{l=0}^{n-m}\frac{\binom{l+m}{l}}{\binom{l+n+1}{l}} S_2(l+n+1,n+1)S_1(n,l+m).
\end{split}
\end{equation*}
\par\bigskip
If is known that
\begin{align*}\tag{7}
x^n\sim (1,t),\quad (x)_n=(1,e^t-1),\quad (see[11]).
\end{align*}
By (1) and (7), we get
\begin{align*}\tag{8}
(x)^n&=x\left(\frac{t}{e^t-1} \right)^nx^{-1}x^n
\\&=x\left(\frac{t}{e^t-1} \right)^nx^{n-1}=xB_{n-1}^{(n)}(x).
\end{align*}
Thus, by (8), we have
\begin{align*}\tag{9}
B_{n-1}^{(n)}(x)=(x-1)_{n-1}, (n \in \mathbb{N}).
\end{align*}
Therefore, by(9), we obtain the following lemma.
\par\bigskip
{\bf Lemma  2 }. For $n \in \mathbb{N}$ , we have
\begin{equation*}
\begin{split}
B_{n}^{(n+1)}(x+1)=(x)_{n}.
\end{split}
\end{equation*}
\par\bigskip
Note that
\begin{align*}\tag{10}
\sum_{l=0}^{\infty}\left(\frac{e^t-1}{t} \right)B_l^{(n+1)}(x)\frac{t^l}{l!}
&=\left(\frac{e^t-1}{t}\right)\left(\frac{t}{e^t-1}\right)^{n+1}e^{xt}
\\&=\left(\frac{t}{e^t-1}\right)^ne^{xt}
\\&=\sum_{l=0}^{\infty}B_l^{(n)}(x)\frac{t^l}{l!}.
\end{align*}
By comparing the coefficients on the both sides of (10), we get
\begin{align*}\tag{11}
\left(\frac{e^t-1}{t}\right)B_l^{(n+1)}(x)=B_l^{(n)}(x),\quad (l \geq 0).
\end{align*}
From (11), we have
\begin{align*}\tag{12}
\left(\frac{e^t-1}{t}\right)B_n^{(n+1)}(x)=B_n^{(n)}(x),\quad (n \geq 0).
\end{align*}
By Lemma 2, (11) and (12), we get
\begin{align*}\tag{13}
B_n^{(n)}(x+1)&=\left(\frac{e^t-1}{t}\right)B_n^{(n+1)}(x+1)
\\&=\left(\frac{e^t-1}{t}\right)(x)_n
\\&= \int_x^{x+1}(u)_ndu.
\end{align*}
From (4) and (13), we have
\begin{align*}\tag{14}
 \int_x^{x+1}(u)_ndu&=\sum_{l=0}^{n}S_1(n,l)\int_{x-1}^{x}u^ldu
\\&=\sum_{l=0}^{n}\frac{S_1(n,l)}{l+1}(x^{l+1}-(x-1)^{l+1}).
\end{align*}
Therefore, by (13) and (14), we obtain the following theorem.
\par\bigskip
{\bf  Theorem 3 }. For $ n\geq1  $ , we have
\begin{equation*}
\begin{split}
B_n^{(n)}(x+1)=\sum_{l=0}^{n}S_1(n,l)\frac{1}{l+1}(x^{l+1}-(x-1)^{l+1}).
\end{split}
\end{equation*}
\par\bigskip
For $a\neq0$, Abel sequence is defined by $A_n (x;a)=x(x-an)^{n-1}$.
\\In [11], we note that $A_n (x;a)\sim(1,te^{at})$.
\\Let us consider the following associated sequences:
\begin{align*}\tag{15}
 &A_{n}(x;a)=x(x-an)^{n-1} \sim (1,te^{at}), \quad a\neq0,
\\ & (\frac{x}{b})_n\sim(1,e^{bt}-1),\quad (b\neq0).
\end{align*}
For $n\geq 1$, by(11), we get
\begin{align*}\tag{16}
 (\frac{x}{b})_n&=x\left(\frac{te^{at}}{e^{bt}-1}\right)^nx^{-1}A_n(x;a)
 \\&=\frac{x}{b^n}\left(\frac{bt}{e^{bt}-1}\right)^ne^{ant}(x-an)^{n-1},
\end{align*}
where
\begin{align*}\tag{17}
 \left(\frac{bt}{e^{bt}-1}\right)^ne^{ant}=\sum_{k=0}^{\infty}b^kB_k^{(n)}\left(\frac{an}{b}\right)\frac{t^k}{k!}.
\end{align*}
From (16) and (17), we have
\begin{align*}\tag{18}
(\frac{x}{b})_n&
=\frac{x}{b^n}\left(\sum_{k=0}^{\infty}b^kB_k^{(n)}\left(\frac{an}{b}\right)\frac{t^k}{k!}\right)(x-an)^{n-1}
\\&=x\sum_{k=0}^{n-1}b^{k-n}B_k^{(n)}\left(\frac{an}{b}\right)\frac{(n-1)_k}{k!}(x-an)^{n-1-k}
\\&=\sum_{k=0}^{n-1}b^{k-n}B_k^{(n)} \left(\frac{an}{b}\right)\binom{n-1}{k}x(x-an)^{n-1-k}
\\&=\sum_{k=0}^{n-1}\binom{n-1}{k}b^{k-n}B_k^{(n)} \left(\frac{an}{b}\right)A_{n-k}(x-ak;a).
\end{align*}
On the other hand, by (16), we get
\begin{align*}\tag{19}
(\frac{x}{b})_n& =\left(\frac{1}{b} \right)^n x \left(\frac{bt}{e^{bt}-1}\right)^n e^{ant}(x-an)^{n-1}
\\&=\left(\frac{1}{b}\right)^n x \left(\frac{bt}{e^{bt}-1}\right)^n x^{n-1}
\\&=\left(\frac{1}{b}\right)^n x \sum_{l=0}^{n-1}\frac{B_l^{(n)}}{l!}b^l(n-1)_lx^{n-l-1}
\\&=\frac{x}{b}\sum_{l=0}^{n-1}\binom{n-1}{l}\left(\frac{1}{b}\right)^{n-l-1}x^{n-l-1}B_l^{(n)}
\\&=\frac{x}{b}B_{n-1}^{(n)}(\frac{x}{b}).
\end{align*}
Therefore, by (18) and (19), we obtain the following theorem.
\par\bigskip
{\bf  Theorem 4 }. For $ n\geq1  $ , we have
\begin{equation*}
\begin{split}
(\frac{x}{b})_n&=\sum_{k=0}^{n-1}\binom{n-1}{k}b^{k-n}B_k^{(n)} \left(\frac{an}{b}\right)A_{n-k}(x-ak;a)
\\&=\frac{x}{b}B_{n-1}^{(n)}(\frac{x}{b}).
\end{split}
\end{equation*}
\par\bigskip
Moreover,
$$xB_{n-1}^{(n)} (\frac{x}{b})=\sum_{k=0}^{n-1}\binom{n-1}{k}b^{k-n+1}B_k^{(n)}\left(\frac{an}{b}\right)A_{n-k}(x-ak;a).$$
{\bf Remark }. For $b=1$, $n\geq 1$, we have
$$ (x)_n=xB_{n-1}^{(n)}(x)=\sum_{k=0}^{n-1}\binom{n-1}{k}B_k^{(n)}(an)A_{n-k}(x-ak;a).$$
Let $\phi_n(x)=\sum_{k=0}^{n}S_2(n,k)x^k$ be exponential polynomial.
\\Then, we note that
\begin{align*}\tag{20}
 \phi_n(x)\sim (1,\log(1+t)), \quad x^n\sim (1,t).
\end{align*}
It is well known that
\begin{align*}\tag{21}
 \left(\frac{\log(1+t)}{t}\right)^n=n\sum_{k=0}^{\infty}\frac{B_k^{(n+k)}}{n+k}\frac{t^k}{k!},(see[5,6]).
\end{align*}
By (1) and (20), we get
\begin{align*}\tag{22}
x^n& =x\left(\frac{\log(1+t)}{t}\right)^nx^{-1}\phi_n(x)
\\&=x\{n\sum_{k=0}^{\infty}\frac{B_k^{(n+k)}}{k+n}\frac{t^k}{k!} \}x^{-1}\phi_n(x)
\\&=n\sum_{k=0}^{n-1}\sum_{l=k+1}^{n}\frac{\binom{l-1}{k}}{n+k}B_k^{(n+k)}S_2(n,l)x^{l-k}
\\&=n\sum_{k=0}^{n-1}\sum_{m=1}^{n-k}\frac{\binom{k+m-1}{k}}{n+k}B_k^{(n+k)}S_2(n,k+m)x^{m}
\\&=n\sum_{m=1}^{n}\{\sum_{k=0}^{n-m}\frac{\binom{k+m-1}{k}}{n+k}B_k^{(n+k)}S_2(n,k+m)\}x^{m}.
\end{align*}
Thus, by (22), we obtain the following theorem.
\par\bigskip
{\bf  Theorem 5 }. For $ n\geq 1  $ with $1\leq m\leq n$, we have
\begin{equation*}
\begin{split}
\sum_{k=0}^{n-m}\frac{n\binom{k+m-1}{k}}{n+k}B_k^{(n+k)}S_2(n,k+m)=\delta_{m,n}.
\end{split}
\end{equation*}
\par\bigskip
Let $M_n(x)=\sum_{k=0}^{n}\binom{n}{k}(n-1)_{n-k}2^k(x)_k$ be Mittag-Leffler polynomials with $M_n(x)\sim(1,\frac{e^t-1}{e^t+1})$.
Then, let us consider the following associated sequence:
\begin{align*}\tag{23}
M_n(x)\sim(1,\frac{e^t-1}{e^t+1}), \quad (x)_n=(1,e^t-1).
\end{align*}
For $n\geq 1$, by (1) and (23), we get
\begin{align*}\tag{24}
(x)_n& =x\left(\frac{1}{e^t+1}\right)^nx^{-1}M_n(x)
\\&=\sum_{k=0}^{n}\binom{n}{k}(n-1)_{n-k}2^kx \left(\frac{1}{e^t+1}\right)^n(x-1)_{k-1}
\\&=\sum_{k=0}^{n}\sum_{l=0}^{k-1}\binom{n}{k}(n-1)_{n-k}2^kS_1(k-1,l)\frac{1}{2^n}x
\left(\frac{2}{e^t+1}\right)^n(x-1)^{l}
\\&=\sum_{k=0}^{n}\sum_{l=0}^{k-1}\binom{n}{k}(n-1)_{n-k}2^{k-n}S_1(k-1,l)x
E_l^{(n)}(x-1).
\end{align*}
Thus, by (24), we obtain the following proposition.
\par\bigskip
{\bf  Proposition 6 }. For $ n\geq 1  $ , we have
\begin{equation*}
\begin{split}
(x)_n=\sum_{k=0}^{n}\sum_{l=0}^{k-1}\binom{n}{k}(n-1)_{n-k}2^{k-n}S_1(k-1,l)x
E_l^{(n)}(x-1).
\end{split}
\end{equation*}
\par\bigskip
For $n\geq 1$, by (1) and (23), we get
\begin{align*}\tag{25}
M_n(x)&=x(e^t+1)^nx^{-1}(x)_n=x(e^t+1)^n(x-1)_{n-1}
\\&=x\sum_{k=0}^{n}\binom{n}{k}e^{kt}(x-1)_{n-1}=x\sum_{k=0}^{n}\binom{n}{k}(x+k-1)_{n-1}.
\end{align*}
The equation (25) is different from the expression
$$M_n(x)=\sum_{k=0}^{n}\binom{n}{k}(n-1)_{n-k}2^k(x)_k.$$
Therefore, by (25), we obtain the following corollary.
\par\bigskip
{\bf  Corollary 7 }. For $ n\geq 1  $ , we have
\begin{equation*}
\begin{split}
M_n(x)=x\sum_{k=0}^{n}\binom{n}{k}(x+k-1)_{n-1}.
\end{split}
\end{equation*}
\par\bigskip
Let $L_n^{(\alpha)}(x)$ be the Laguerre polynomials of order $\alpha (\in \mathbb{R})$. Then we note that $L_n^{(\alpha)}(x)\sim ((1-t)^{-\alpha-1},\frac{t}{t-1})$. Especially, $L_n(x)\sim (1,\frac{t}{t-1})$. By the definition of associated sequences, we see that
\begin{align*}\tag{26}
<\left( \frac{t}{t-1}\right )^n | \quad L_n{(x)}>=n!\delta _{n,k},\quad (n,k \geq 0).
\end{align*}
From (26), we have
\begin{align*}\tag{27}
<\left( \frac{t}{t+1}\right )^n | \quad L_n {(-x)}>=n!\delta _{n,k}.
\end{align*}
Thus, by (27), we get
\begin{align*}\tag{28}
 L_n {(-x)} \sim (1,\frac{t}{t+1}).
\end{align*}
As it is shown in Roman [11], one can find an explicit expression for
$L_n(x)$ by using the transfer formula
\begin{align*}\tag{29}
 L_n(-x)=\sum_{k=1}^{n}\binom{n-1}{k-1}\frac{n!}{k!}x^k,\quad (n\geq 1), (see[11]).
\end{align*}
It is well known that
\begin{align*}\tag{30}
 \frac{t}{(1+t)\log(1+t)}=\sum_{k=0}^{\infty}B_k^{(k)}\frac{t^k}{k!}, (see[5,6]).
\end{align*}
Thus, by (29), we get
\begin{align*}\tag{31}
\left( \frac{t}{(1+t)\log(1+t)}\right)^n=\sum_{k=0}^{\infty}
\left(\sum_{l_1+...+l_n=k}^{} \binom{k}{l_1,...,l_n}B_{l_1}^{(l_1)}...B_{l_n}^{(l_n)}\right)\frac{t^k}{k!}.
\end{align*}
From (1),(20) and (28), we have
\begin{align*}\tag{32}
\phi_n(x)&=x\left( \frac{t}{(1+t)\log(1+t)}\right)^nx^{-1}L_n(-x)
\\&=x\left( \frac{t}{(1+t)\log(1+t)}\right)^nx^{-1}
\sum_{m=1}^{n}\binom{n-1}{m-1}\frac{n!}{m!}x^m.
\end{align*}
By (30) and (31), we get
\begin{align*}\tag{33}
&\phi_n(x) \\&=\sum_{m=1}^{n}\binom{n-1}{m-1}\frac{n!}{m!}x
\{\sum_{k=0}^{m-1}\sum_{l_1+...+l_n=k}^{}
\binom{k}{l_1,...,l_n}B_{l_1}^{(l_1)}...B_{l_n}^{(l_n)}
\}\frac{t^k}{k!}x^{m-1}
\\&=\sum_{m=1}^{n}\sum_{k=0}^{m-1}\sum_{l_1+...+l_n=k}^{}\binom{n-1}{m-1}\binom{m-1}{k}\frac{n!}{m!}
\binom{k}{l_1,...,l_n}B_{l_1}^{(l_1)}...B_{l_n}^{(l_n)}x^{m-k}
\\&=\sum_{m=1}^{n}\{\sum_{l=1}^{m}\sum_{l_1+...+l_n=m-l}^{}\binom{n-1}{m-1}\binom{m-1}{l-1}\frac{n!}{m!}
\binom{m-l}{l_1,...,l_n}B_{l_1}^{(l_1)}...B_{l_n}^{(l_n)}\}x^{l}.
\end{align*}
From (20), we have
\begin{align*}\tag{34}
\phi_n(x)=\sum_{k=0}^{n}S_2(n,k)x^k=\sum_{k=1}^{n}S_2(n,k)x^k, \quad (n\geq 1).
\end{align*}
Therefore, by (32) and (33), we obtain the following theorem.
\par\bigskip
{\bf  Theorem 8 }. For $ n\geq 1  $ with $1\leq l\leq n$, we have
\begin{equation*}
\begin{split}
&S_2(n,l)
\\&=\sum_{l\leq m\leq
n}^{}\sum_{l_1+...+l_n=m-l}^{}\binom{n-1}{m-1}\binom{m-1}{l-1}\frac{n!}{m!}
\binom{m-l}{l_1,...,l_n}B_{l_1}^{(l_1)}...B_{l_n}^{(l_n)}.
\end{split}
\end{equation*}
\par\bigskip
It is well known in [5,6] that
\begin{align*}\tag{35}
\frac{(e^t-1)^n}{e^{tx}t^n}=(n!)^2\sum_{k=0}^{\infty}\left(\sum_{j=0}^{k}(-1)^{k-j}
\frac{\binom{k}{j}S_2(j+n,n)}{\binom{j+n}{j}}x^{k-j} \right)\frac{t^k}{k!}.
\end{align*}
From (35), we have
\begin{align*}\tag{36}
\left(\frac{e^{bt}-1}{te^{at}}\right)^n=(n!)^2b^n\sum_{k=0}^{\infty}\left(\sum_{j=0}^k(-1)^{k-j}
\frac{\binom{k}{j}S_2(j+n,n)}{\binom{j+n}{j}}(an)^{k-j}b^j
\right)\frac{t^k}{k!},
\end{align*}
where $a,b \neq 0$.
\\By (1) and (15), we get
\begin{align*}\tag{37}
A_n(x;a)&= x\left(\frac{e^{bt}-1}{te^{at}}\right)^nx^{-1}(\frac{x}{b})_n
\\&= x\left(\frac{e^{bt}-1}{te^{at}}\right)^n\frac{1}{b}(\frac{x}{b}-1)_{n-1}
\\&=(n!)^2b^{n-1}x\sum_{k=0}^{n-1}\left(\sum_{j=0}^{k}(-an)^{k-j}b^j
\frac{\binom{k}{j}S_2(j+n,n)}{\binom{j+n}{j}}\right)\frac{t^k}{k!}(\frac{x}{b}-1)_{n-1},
\end{align*}
where
\begin{align*}\tag{38}
t^k
(\frac{x}{b}-1)_{n-1}&=\sum_{l=0}^{n-1}S_1(n-1,l)t^k(\frac{x}{b}-1)^{l}
 \\&=\sum_{l=k}^{n-1}S_1(n-1,l)(\frac{1}{b})^k(l)_k (\frac{x}{b}-1)^{l-k}.
\end{align*}
From (36) and (37), we have
\begin{align*}\tag{39}
&A_n(x;a)=x(x-an)^{n-1}
\\&=(n!)^2b^{n-1}\sum_{k=0}^{n-1}\left(\sum_{j=0}^{k}\sum_{l=k}^{n-1}\left(-\frac{an}{b}\right)^{k-j}
\frac{\binom{k}{j}\binom{l}{k}S_2(j+n,n)S_1(n-1,l)}{\binom{j+n}{j}}\right)x(\frac{x}{b}-1)^{l-k}.
\end{align*}
Therefore, by (39), we obtain the following lemma.
\par\bigskip
{\bf  Lemma 9 }. For $ n\geq 1  $, we have
\begin{equation*}
\begin{split}
&A_n(x;a)
\\&=(n!)^2b^{n-1}\sum_{k=0}^{n-1}\left(\sum_{j=0}^{k}\sum_{l=k}^{n-1}\left(-\frac{an}{b}\right)^{k-j}
\frac{\binom{k}{j}\binom{l}{k}S_2(j+n,n)S_1(n-1,l)}{\binom{j+n}{j}}\right)x(\frac{x}{b}-1)^{l-k}.
\end{split}
\end{equation*}
\par\bigskip
{\bf Remark.} Let $b=1$. Then we have
\begin{align*}
&A_n(x;a)=x(x-an)^{n-1}
\\&=(n!)^2\sum_{k=0}^{n-1}\left(\sum_{j=0}^{k}\sum_{l=k}^{n-1}(-an)^{k-j}
\frac{\binom{k}{j}\binom{l}{k}S_2(j+n,n)S_1(n-1,l)}{\binom{j+n}{j}}\right)x(x-1)^{l-k}.
\end{align*}
It is well known in[5,6] that
\begin{align*}\tag{40}
 \frac{(1+t)^{x-1}t^n}{(\log(1+t) )^n}=\sum_{k=0}^{\infty}B_k^{(k-n+1)}(x)\frac{t^k}{k!}.
\end{align*}
From (39), we have
\begin{align*}\tag{41}
 \left(\frac{t(1+t)^{a}}{\log(1+t)}\right)^n=\sum_{k=0}^{\infty}B_k^{(k-n+1)}(an+1)\frac{t^k}{k!}.
\end{align*}
Let us consider the following associated sequences:
\begin{align*}\tag{42}
 S_n {(x)} \sim (1,t(t+1)^a), (a\neq0) ,\quad  x^n \sim (1,t).
\end{align*}
Then, for $n\geq1$, by (1) and (41), we get
\begin{align*}\tag{43}
S_n(x)&=x\left( \frac{t}{t(1+t)^a }\right)^nx^{-1}x^n
\\&=x(1+t)^{-an}x^{n-1}
\\&=x\sum_{k=0}^{n-1}\binom{-an}{k}(n-1)_kx^{n-1-k}
\\&=\sum_{k=1}^{n}\binom{-an}{n-k}(n-1)_{n-k}x^{k}.
\end{align*}
Therefore, by (42), we obtain the following proposition.
\par\bigskip
{\bf  Proposition 10 }. For $ n\geq 1  $, let $S_n(x) \sim (1,t(1+t)^a)$, $(a\neq0)$.
\\Then, we have
\begin{equation*}
\begin{split}
S_n(x)=\sum_{k=1}^{n}\binom{-an}{n-k}(n-1)_{n-k}x^{k}.
\end{split}
\end{equation*}
\par\bigskip
By (1), (20) and (41), we get
\begin{align*}\tag{44}
\phi_n(x)&=x\left( \frac{t(1+t)^a}{\log(1+t)}\right)^nx^{-1}S_n(x)
\\&=x\left( \frac{t(1+t)^a}{\log(1+t)}\right)^n
\sum_{l=1}^{n}\binom{-an}{n-l}(n-1)_{n-l}x^{l-1}.
\end{align*}
From (40) and (44), we get
\begin{align*}\tag{45}
\phi_n(x)&=\sum_{l=1}^{n}\binom{-an}{n-l}(n-1)_{n-l}x \left(\frac{t(1+t)^a}{\log(1+t)}\right)^nx^{l-1}
\\&=\sum_{l=1}^{n}\binom{-an}{n-l}(n-1)_{n-l}\sum_{k=0}^{l-1}B_k^{(k-n+1)}(an+1)\frac{1}{k!}xt^kx^{l-1}
\\&=\sum_{l=1}^{n}\sum_{k=0}^{l-1}\binom{-an}{n-l}(n-1)_{n-l}\binom{l-1}{k}B_k^{(k-n+1)}(an+1)x^{l-k}
\\&=\sum_{l=1}^{n}\sum_{m=1}^{l}\binom{-an}{n-l}(n-1)_{n-l}\binom{l-1}{m-1}B_{l-m}^{(l-n-m+1)}(an+1)x^{m}
\\&=\sum_{m=1}^{n}\{\sum_{l=m}^{n}\binom{-an}{n-l}(n-1)_{n-l}\binom{l-1}{m-1}B_{l-m}^{(l-n-m+1)}(an+1)\}x^{m}.
\end{align*}
Therefore, by (20) and (44), we obtain the following theorem.
\par\bigskip
{\bf  Theorem 11 }. For $ n\geq 1  $ and $m\geq 0$, we have
\begin{equation*}
\begin{split}
S_2(n,m)=\sum_{m\leq l\leq
n}^{}\binom{-an}{n-l}(n-1)_{n-l}\binom{l-1}{m-1}B_{l-m}^{(l-m-n+1)}(an+1).
\end{split}
\end{equation*}
\par\bigskip
Let us consider the following associated sequences:
\begin{align*}\tag{46}
S_n^*(x)\sim(1, \frac{t(e^t+1)}{2}) , \quad x^n\sim(1,t)
\end{align*}
Then, by (1) and (46), we get
\begin{align*}\tag{47}
S_n^*(x)=x \left(\frac{2}{e^t+1}\right)^nx^{n-1}=xE_{n-1}^{(n)}(x).
\end{align*}
For $n\geq1$, by (1), (15) and (46), we get
\begin{align*}\tag{48}
xE_{n-1}^{(n)}(x)&=x \left(\frac{2}{e^t+1}\right)^ne^{ant}(x-an)^{n-1}
\\&=x\sum_{k=0}^{n-1}\frac{E_k^{(n)}(an)}{k!}t^k(x-an)^{n-1}
\\&=\sum_{k=0}^{n-1}\binom{n-1}{k} E_k^{(n)}(an)x(x-an)^{n-1-k}.
\end{align*}
Therefore, by (47) , we obtain the following theorem.
\par\bigskip
{\bf  Theorem 12 }. For $ n\geq 1 $ , we have
\begin{equation*}
\begin{split}
xE_{n-1}^{(n)}(x)=\sum_{k=0}^{n}\binom{n-1}{k}E_k^{(n)}(an)A_{n-k}(x-ak;a).
\end{split}
\end{equation*}
\par\bigskip
ACKNOWLEDGEMENTS. This paper was supported by Basic Science Research Program through the National  Research Foundation of Korea(NRF)
funded by the Ministry of Education, Science and Technology 2012R1A1A2003786. The second author was supported in part by the Research Grant of Kwangwoon University in 2013.

\par\noindent
{\mpn { \bpn {\small Dae San {\sc Kim} \mpn Department of Mathematics,
\pn Sogang University, Seoul 121-742, S. Korea \pn {\it E-mail:}\
{\sf dskim@sogang.ac.kr} }
 \mpn { \bpn {\small Taekyun {\sc KIM} \mpn
 Department of Mathematics,\pn
Kwangwoon University, Seoul 139-701, S.Korea
  \pn {\it E-mail:}\ {\sf tkkim@kw.ac.kr} }
\mpn { \bpn {\small Seog-Hoon {\sc Rim} \mpn Department of
Mathematics Education, \pn Kyungpook National University, Taegu
702-701, S. Korea \pn {\it E-mail:}\ {\sf shrim@knu.ac.kr} }

\

\end{document}